\documentclass{article}
\usepackage{dsfont}
\title{$MPLS$ $=$ Mutually Projective Latin Squares}
\author{Leendert Bleijenga}
\date{6 March 2012}
\begin{document}
\maketitle

The author is born on 14 March 1940.\footnote{email: leen.bleijenga@gmail.com}

\begin{abstract}
We will see that every finite projective plane of order $\kappa \geq 2$ gives rise to a complete set of $(\kappa - 1)$ $MPLS$ (= mutually projective latin squares)\footnote{\`{A} bas les MOLS = mutually orthogonal latin squares.} of order $\kappa$ and by reversing the process we can construct a finite projective plane of order $\kappa$ when a complete set of $(\kappa - 1)$ $MPLS$ of order $\kappa$ is given.
\end{abstract}

\setlength{\parskip}{5mm}

Subjects: Combinatorics

\newpage

\section{Finite Geometries}

A projective plane is a geometry. When is a geometry a projective plane? We first give a definition of a geometry.

\newtheorem{defi}{Definition}
\begin{defi}
A pair $(\mathcal{P}, \mathcal{L})$ where $\mathcal{P}$ and $\mathcal{L}$ are sets, and the elements of $\mathcal{L}$ are subsets of $\mathcal{P}$, is called a \textbf{geometry} with point set $\mathcal{P}$ and line set $\mathcal{L}$ (we assume $\mathcal{P} \cap \mathcal{L} =\emptyset$) iff

\begin{enumerate}
  \item Every set $\{p, q\}$ consisting of two different points in $\mathcal{P}$ is a subset of exactly one line, let us say $pq = qp = \ell \in \mathcal{L}$;
  \item Every line $m \in \mathcal{L}$ contains at least two different points of $\mathcal{P}$.
\end{enumerate}
\end{defi}

\newtheorem{theodef}[defi]{Theorem and definition}
\begin{theodef}
Let $(\mathcal{P}, \mathcal{L})$ be a geometry and let $\mathcal{Q} \subset \mathcal{P}$. Define the set $\mathcal{M} = \{\mathcal{Q} \cap \ell | \ell \in \mathcal{L} \  \& \  |\mathcal{Q} \cap \ell| \geq 2 \}$. Then the pair $(\mathcal{Q}, \mathcal{M})$ is a geometry and is called a \textbf{subgeometry} of $(\mathcal{P}, \mathcal{L})$.
\end{theodef}

\begin{defi}
A subgeometry $(\mathcal{Q}, \mathcal{M})$ of the geometry $(\mathcal{P}, \mathcal{L})$ is called a \textbf{subspace} of $(\mathcal{P}, \mathcal{L})$ if $\mathcal{Q} \subset \mathcal{P}$ and $\mathcal{M} \subset \mathcal{L}$.
\end{defi}

\begin{defi}
Examples of subspaces of the geometry $(\mathcal{P}, \mathcal{L})$ are $(\emptyset, \emptyset)$ (this is called the \textbf{empty space} or the \textbf{vacuum}), $(\{p\}, \emptyset)$ where $p \in \mathcal{P}$ (these are called the \textbf{point spaces}), $(\ell, \{\ell\})$ where $\ell \in \mathcal{L}$ (these are called the \textbf{line spaces}) and $(\mathcal{P}, \mathcal{L})$ (this is called the \textbf{universe}).
\end{defi}

\newtheorem{theo}[defi]{Theorem}
\begin{theo}
Let $(\mathcal{P}, \mathcal{L})$ be a geometry with point $p$ and line $\ell$ such that point $p$ is not incident with line $\ell$, $p \not \in \ell$ and let $a, b$ be points incident with line $\ell$, $a, b \in \ell$ and $a \not = b$. Then $pa \not = pb$.
\end{theo}

\textbf{Proof:} If $pa = pb$ then $p, a, b \in pa = pb$. Then $p \in ab = \ell$. A contradiction.

\begin{defi}
A geometry $(\mathcal{P}, \mathcal{L})$ is called a \textbf{plane} if the maximal subspaces of $(\mathcal{P}, \mathcal{L})$ are the line spaces.
\end{defi}

\begin{defi}
A geometry $(\mathcal{P}, \mathcal{L})$ is called \textbf{finite} if $|\mathcal{P}| \in \mathds{N}$.
\end{defi}

In a finite geometry we define $|\mathcal{P}| = v$ ($v$ for varieties) and $|\mathcal{L}| = b$ ($b$ for blocks).

\begin{defi}
A finite geometry $(\mathcal{P}, \mathcal{L})$ is called \textbf{regular} if each point of $\mathcal{P}$ lies on the same number $r \in \mathds{N}$ of lines ($r$ for replications).
\end{defi}

\begin{defi}
A finite geometry $(\mathcal{P}, \mathcal{L})$ is called \textbf{uniform}
if each line of $\mathcal{L}$ has the same number $k \in \mathds{N}$ of points ($k$ for plots).
\end{defi}

\newtheorem{theo1}[defi]{Theorem of de Bruijn \& Erd\"os [1948]}
\begin{theo1}
Let $(\mathcal{P}, \mathcal{L})$ be a finite geometry and let $|\mathcal{L}| \geq 2$. Then $v \leq b$. If $v = b$ then any two different lines of $\mathcal{L}$ intersect and $(a)$ the geometry consists of a pencil and a transversal not through the top of the pencil or $(b)$ the geometry is a projective plane. [$(c)$ If any two different lines of $\mathcal{L}$ intersect then $v = b$.]
\end{theo1}

\textbf{Proof:} (Footnote 1 of [\textbf{1948}] states: This [i.e. $v \leq b$] was also proved by G. SZEKERES but his proof was more complicated.)
\newline We follow the combinatorial proof of [\textbf{1948}]:
\newline Let $\mathcal{P} = \{p_{1}, p_{2}, \cdots, p_{v}\}$ and let $\mathcal{L} = \{\ell_{1}, \ell_{2}, \cdots, \ell_{b}\}$, $b \geq 2$. So $v \geq 3$. Let $r_{i}$ denote the number of lines passing through point $p_{i}, 1 \leq i \leq v$, and let $k_{j}$ denote the number of points lying on line $\ell_{j}, 1 \leq j \leq b$. Now the number of point-line incidences is the same as the number of line-point incidences, so
  \begin{equation}
  \sum_{i = 1}^{v} r_{i} = \sum_{j = 1}^{b} k_{j}
  \end{equation}
If point $p_{i}$ is not incident with line $\ell_{j}$ then
  \begin{equation}
  r_{i} \geq k_{j}
  \end{equation}
because $p_{i}$ can be connected not only with the $k_{j}$ points of line $\ell_{j}$, these $k_{j}$ connection lines  are different, but through $p_{i}$ there might be lines which are parallel to $\ell_{j}$.
Let $r_{v} = w$ be the smallest number such that $r_{i} \geq r_{v}, i = 1, 2, \cdots, v$, and let $\ell_{1}, \ell_{2}, \cdots, \ell_{w}, 2 \leq w \leq b$, be the $w$ lines through point $p_{v}$. We can choose $w$ points on these $w$ lines and we can call them $p_{1}, p_{2}, \cdots, p_{w}$ all different from point $p_{v}$. We have
  \begin{equation}
  r_{1} \geq k_{2}, r_{2} \geq k_{3}, \cdots, r_{w - 1} \geq k_{w}, r_{w} \geq k_{1}, r_{i} \geq r_{v} = w \geq k_{j}, i, j > w
  \end{equation}
for otherwise there are more than $w$ lines through point $p_{v}$.
Now if $b < v$ then all the terms on the right-hand side of $(1)$ can be majorized qua $\geq$ by $b$ terms of the left-hand side of $(1)$ by $(3)$, so $(1)$ would become an inequality by taking the remaining $v - b$ terms of the left-hand side of $(1)$ into account. From this contradiction it follows $v \leq b$.

In case of $v = b$ all inequalities of $(3)$ becomes equalities:
  \begin{equation}
  r_{1} = k_{2}, r_{2} = k_{3}, \cdots, r_{w - 1} = k_{w}, r_{w} = k_{1}, r_{i} = r_{v} = w = k_{j}, i, j > w
  \end{equation}
Without loss of generality we can assume $r_{1} \leq r_{2} \leq \cdots \leq r_{w}$.
\begin{itemize}
  \item[(a)] Let $r_{1} < r_{w}$, which implies $k_{2} < k_{1}$. Thus $k_{1} = r_{w} 
  > r_{1} \geq r_{v} = w$. For all points $p_{i}$ and lines $\ell_{j}, i,j > w$, we have $r_{i} = w = k_{j}$. From this it follows that the lines $\ell_{2}, \ell_{3}, \cdots, \ell_{w}$ all contains exactly two different points, otherwise there exists a point $p_{i}, i > w$, on $\ell_{w}$, say, such that $p_{i} \neq p_{v}, p_{w}$ but then $r_{i} > w$, by connecting point $p_{i}$ with the points on line $\ell_{1}$. A contradiction. On the connection line $\ell_{x}$ of the points $p_{1}$ and $p_{2}$ there are $w$ points which must be the points $p_{1}, p_{2}, \cdots, p_{w}$. Let the point $p_{w + 1}$ lie on line $\ell_{1}$ but different from the points $p_{1}$ and $p_{v}$. On the connection line of the points $p_{w}$ and $p_{w+ 1}$ there are $w$ points too which must be $p_{2}, p_{3}, \cdots, p_{w}$ and $p_{w + 1}$, but if $w \geq 3$ then we have $p_{1}p_{2} = p_{w}p_{w + 1} = p_{1}p_{w + 1} = \ell_{1}$ so all points are collinear. A contradiction. Thus $w = 2$. In this case there are $v - 1$ collinear points, lets us say $p_{1}, p_{2}, \cdots, p_{v - 1},$ lying on the transversal, lets us say $\ell_{b}$, and a point outside the transversal, lets us say $p_{v}$, which is the top of a pencil consisting of $v - 1$ connection lines, connecting the top with the points of the transversal.
  \item[(b)] Let $r_{1} = r_{w} \geq w$. Then $r_{1} = r_{2} = \cdots = r_{w} \geq w$ and $k_{2} = k_{3} = \cdots = k_{1} = r_{w} = r_{1} \geq w$. For all points $p_{i}$ and lines $\ell_{j}, i,j > w$, we have $r_{i} = w = k_{j}$. Suppose that point $p_{w + 1}$ lies on line $\ell_{2}$ and that $p_{v} \neq p_{w + 1} \neq p_{2}$ and that $k_{1} > w$. Then $r_{w + 1} = w$ but if we connect point $p_{w + 1}$ with all the points of $\ell_{1}$ then we have $r_{w + 1} > w$. A contradiction. Thus $k_{1} = w$ and the geometry is regular with replication number $r$ and uniform with plot number $k$ and $r = k = w$. We show that any two lines intersect. Let $\ell_{x}$ and $\ell_{y}$ be two different lines that are parallel, with point $p_{x}$ on line $\ell_{x}$ and point $p_{y}$ on line $\ell_{y}$. Then $k_{x} = k_{y} = w$. But if we connect point $p_{x}$ with the $w$ points of line $\ell_{y}$, we see that $r_{x} \geq w + 1$. A contradiction. Thus two different lines intersect and the geometry is a projective plane for $w > 2$.
  \item[(c)] [The proof of the last part of the theorem is left to the reader. (Consider the case where the maximum number of independent points is three and the case where there are at least four independent points.)]
\end{itemize}

\newtheorem{theo2}[defi]{Theorem of Bleijenga [1993]}
\begin{theo2}
Let $(\mathcal{P}, \mathcal{L})$ be a finite geometry, and let $|\mathcal{L}| \geq 2$ then there exists an injection $\sigma: \mathcal{P} \rightarrow \mathcal{L}$ such that $p \in \sigma p$ for all $p \in \mathcal{P}$.
\end{theo2}

\textbf{Proof:} See proof of R.H. Jeurissen [1995].

\begin{defi}
A set of points of a geometry are called \textbf{collinear} if these points lie on the same line.
\end{defi}

\begin{defi}
A set of points of a geometry are called \textbf{independent} if no three points lie on the same line.
\end{defi}

\begin{theo}
  Let $(\mathcal{P}, \mathcal{L})$ be a finite geometry and let $|\mathcal{L}| \geq 2$. Let $\ell_{1}$ and $\ell_{2}$ be two different lines which either intersect in a point $p$ or are parallel (= disjoint). Let $p_{1}$ and $p_{2}$ be two different points on line $\ell_{1}$ different from point $p$ and let $p_{3}$ and $p_{4}$ be two different points on line $\ell_{2}$ different from point $p$. Then the points $p_{1}$, $p_{2}$, $p_{3}$ and $p_{4}$ are four independent points.
\end{theo}

We now are ready to give two definitions of a projective plane.

\newtheorem{defi1}[defi]{First definition}
\begin{defi1}
A finite geometry $(\mathcal{P}, \mathcal{L})$ is called a \textbf{projective plane} if
\begin{enumerate}
  \item it is \textbf{regular} and \textbf{uniform};
  \item $r = k$;
  \item there are at least four independent points.
\end{enumerate}
\end{defi1}

\newtheorem{defi2}[defi]{Second definition}
\begin{defi2}
A finite geometry $(\mathcal{P}, \mathcal{L})$ is called a \textbf{projective plane} if
\begin{enumerate}
  \item two different lines intersect;
  \item there are at least four independent points.
\end{enumerate}
\end{defi2}

Bear in mind that the empty space $(r = k = 0)$, the point spaces $(r = k = 0)$ and the triangle spaces with three points and three lines $(r = k = 2)$ are regular and uniform but are not allowed as projective planes by both the First and the Second definition $(15)$ and $(16)$.

\begin{defi}
The order $\kappa$(kappa) of a projective plane is equal to the number of points $k$ on a line minus one, so $k = \kappa + 1$.
\end{defi}

\begin{theo}
A projective plane of order $\kappa$ has $v = b = \kappa^{2} +\kappa + 1$ points and lines.
\end{theo}

\section{Canonical Incidence Matrices}

The incidence matrix $M = [m_{ij}]$ of order $\kappa^{2} + \kappa + 1$ of a projective plane $(\mathcal{P}, \mathcal{L})$ of order $\kappa \geq 2$, where $m_{ij} = 1$ when line $\ell_{i}$ is incident with point $p_{j}$, and $m_{ij} = 0$ otherwise, $1 \leq i, j \leq \kappa^{2} + \kappa + 1$, can be put in an easy to handle form after suitable permutations of the rows and permutations of the columns. We call the matrix after permutation also $M$. Firstly we permute the columns of matrix $M$ in such a way that the first $\kappa + 1$ columns have solely ones in the first row and the last $\kappa^{2}$ columns of $M$ have solely zeroes in the first row. Next we permute the rows of matrix $M$ (the first row remains unchanged however) in such a way that the first $\kappa + 1$ rows have solely ones in the first column and the last $\kappa^{2}$ rows of $M$ have solely zeroes in the first column. Secondly we permute the columns of $M$, leaving the first $\kappa + 1$ columns intact, such that we have in the second row solely ones in the columns $\kappa + 2$ through $2 \kappa + 1$ and solely zeroes in the last $(\kappa^{2} - \kappa)$ columns. x-ly, $2 \leq x \leq \kappa + 1$ we permute the columns of matrix $M$, leaving the first $(x - 1) \kappa + 1$ columns intact, such that we have in the $x$-th row solely ones in the columns $(x - 1)\kappa + 2$ through $x \kappa + 1$ and solely zeroes in the last $\kappa^{2} -(x - 1)\kappa$ columns. Dually we permute the rows of $M$, leaving the first $\kappa + 1$ rows intact, such that we have in the second column solely ones in the rows $\kappa + 2$ through $2 \kappa + 1$ and solely zeroes in the last $(\kappa^{2} - \kappa)$ rows. x-ly, $2 \leq x \leq \kappa + 1$ we permute the rows of matrix $M$, leaving the first $(x - 1) \kappa + 1$ rows intact, such that we have in the $x$-th column solely ones in the rows $(x - 1)\kappa + 2$ through $x \kappa + 1$ and solely zeroes in the last $\kappa^{2} - (x - 1)\kappa$ rows. We partition the matrix $M = [A_{ij}]$ of order $\kappa + 1$ as follows: submatrix $A_{11}$ of order $\kappa + 1$ contains solely ones in the first row and the first column and zeroes elsewhere. The submatrices $A_{1r}$, $2 \leq r \leq \kappa + 1$, of order $(\kappa + 1) \times \kappa$, have solely ones in row $r$ and zeroes elsewhere. The submatrices $A_{s1}$, $2 \leq s \leq \kappa + 1$, of order $\kappa \times (\kappa + 1)$, have solely ones in column $s$ and zeroes elsewhere. The remaining submatrices $A_{ij}$, $2 \leq i, j \leq \kappa + 1$ are of order $\kappa$ which we now investigate further. 

\begin{theo}
The $\kappa^{2}$ submatrices $A_{ij}$, $2 \leq i, j \leq \kappa + 1$ are permutation matrices of order $\kappa$.
\end{theo}

Proof: Every row (= line) and every column (= point) of matrix $M$ contains $\kappa + 1$ ones. This means that every row and every column of the square submatrix $N = [A_{ij}], 2 \leq i, j \leq \kappa + 1$ of $M$ contains $\kappa$ ones. Any of the submatrices $A_{ij}$ has in each row at most one one, otherwise the matrix $A_{1j}$ has two ones in the same columns, but then there are two lines having two points in common, a contradiction. The matrices $A_{ij}$ must have exactly one one in each row to add up to $\kappa$ in the same row (pigeonhole principle). Dually the submatrices $A_{ij}, 2 \leq i, j \leq \kappa + 1$ must have exactly one one in each column. So each submatrix $A_{ij}, 2 \leq i, j \leq \kappa + 1$ is a permutation matrix.

By permuting the columns of the submatrices $A_{2j}$, $2 \leq j \leq \kappa + 1$ we can assume $A_{2j} = I_{\kappa}$, where $I_{\kappa}$ is the identity matrix of order $\kappa$. Dually, by permuting the rows of the submatrices $A_{i2}$, $3 \leq i \leq \kappa + 1$ we can assume $A_{i2} = I_{\kappa}$.

\begin{theo}
Two different submatrices $A_{ij_{1}}$ and $A_{ij_{2}}$, $2 \leq j_{1}, j_{2} \leq \kappa + 1$, $3 \leq i \leq \kappa + 1$ have no one in the same position (= row and column).
\end{theo}

Proof: The two columns of matrix $M$ through the two ones in the same row and the same column position are also going through the two ones in the same row and same column position in the submatrices $A_{2j_{1}} = I_{\kappa}$ and $A_{2j_{2}} = I_{\kappa}$ and the two rows have two ones in common which is a contradiction.

\begin{theo}
Two different submatrices $A_{i_{1}j}$ and $A_{i_{2}j}$, $2 \leq i_{1}, i_{2} \leq \kappa + 1$, $3 \leq j \leq \kappa + 1$ have no one in the same position.
\end{theo}

Proof: Dually the same way as the proof of theorem $20$.

\begin{theo}
  $\sum_{j = 2}^{\kappa + 1} A_{ij} = J, i = 3, 4, \cdots, (\kappa + 1)$, where the elements of $J$ are all one.
\end{theo}

\begin{theo}
  $\sum_{i = 2}^{\kappa + 1} A_{ij} = J, j = 3, 4, \cdots, (\kappa + 1)$, where the elements of $J$ are all one.
\end{theo}

We define the matrix $L_{i} = \sum_{j = 2}^{\kappa + 1} (j - 1) A_{(i + 2)j}, i = 1, 2, \cdots, (\kappa - 1)$.

\begin{theo}
The matrices $L_{i}, 1 \leq i \leq \kappa - 1$ are latin squares.
\end{theo}

Proof: This follows from theorems $22$ and $23$.

\begin{defi}
Two latin squares $L_{1} = [l_{ij}^{(1)}]$ and $L_{2} = [l_{ij}^{(2)}]$ of order $\kappa$ are called \textbf{projective} if they have solely ones in the diagonal and if any row of $L_{1}$ and any row of $L_{2}$ has exactly one element in common which is in the same column of both rows.
\end{defi}

\begin{defi}
A number of latin squares of order $\kappa$ are called $MPLS$ (= \textbf{mutually projective latin squares} if they are pairwise projective.
\end{defi}

\begin{theo}
There are at most $\kappa - 1$ $MPLS$ (= mutually projective latin squares) of order $\kappa$.
\end{theo}

Proof: The elements in the first row and first column are all equal to one. The elements in the first row and second  column must all be different and can be $2, 3, \cdots, \kappa$. So the number is at most $\kappa - 1$.

\begin{defi}
A set of $MPLS$ (= mutually projective latin squares) of order $\kappa$ is \textbf{complete} if the set contains $\kappa - 1$ latin squares.
\end{defi}

\begin{theo}
If a complete set of $\kappa - 1$ $MPLS$ $L_{i}$, $1 \leq i \leq \kappa - 1$ of order $\kappa$ is given then there exists a projective plane of order $\kappa$.
\end{theo}

Proof: The incidence matrix $M$ is easily reconstructed by reversing the entire procedure.

\begin{theo}
In the columns $i$ and $j$, $i \neq j$ of all $\kappa(\kappa - 1)$ rows of a complete set of mutually projective latin squares $L_{1}, L_{2}, \cdots, L_{\kappa - 1}$ of order $\kappa$, we find all possible $\kappa(\kappa - 1)$ variations two by two of $\kappa$ elements $1, 2, \cdots \kappa$.
\end{theo}

Remark: If $L_{1}, L_{2}, \cdots, L_{\kappa - 1}$ is a complete set of mutually projective latin squares (MPLS) of order $\kappa$ then the transposes $L^{T}_{1}, L^{T}_{2}, \cdots, L^{T}_{\kappa - 1}$ need not to be a complete set of MPLS.

The following is an example of a complete set of $4$ $MPLS$ of order $5$:

$\left[ \begin{array}{ccccc}
1 & 2 & 3 & 4 & 5 \\
4 & 1 & 5 & 3 & 2 \\
5 & 3 & 1 & 2 & 4 \\
2 & 5 & 4 & 1 & 3 \\
3 & 4 & 2 & 5 & 1
\end{array} \right], 
\left[ \begin{array}{ccccc}
1 & 3 & 4 & 5 & 2 \\
5 & 1 & 2 & 4 & 3 \\
2 & 4 & 1 & 3 & 5 \\
3 & 2 & 5 & 1 & 4 \\
4 & 5 & 3 & 2 & 1
\end{array} \right], 
\left[ \begin{array}{ccccc}
1 & 4 & 5 & 2 & 3 \\
2 & 1 & 3 & 5 & 4 \\
3 & 5 & 1 & 4 & 2 \\
4 & 3 & 2 & 1 & 5 \\
5 & 2 & 4 & 3 & 1
\end{array} \right]$,

$\left[ \begin{array}{ccccc}
1 & 5 & 2 & 3 & 4 \\
3 & 1 & 4 & 2 & 5 \\
4 & 2 & 1 & 5 & 3 \\
5 & 4 & 3 & 1 & 2 \\
2 & 3 & 5 & 4 & 1
\end{array} \right]$.

\begin{defi}
Let $N = [n_{ij}]$ be a matrix of order $\kappa$. Some elements of $N$ are called \textbf{independent} if no two of them lie in the same row or the same column.
\end{defi}

Remark: The matrix $N$ contains $\kappa$ independent elements.

\begin{defi}
Let $N$ be a latin square of order $\kappa$. A \textbf{transversal} of $N$ is a matrix $T$ of order $\kappa$ with a copy of $\kappa$ different and independent elements of $N$ equal to $1, 2, \cdots, \kappa$ and the other elements equal to zeroes.
\end{defi}

Remark: If we replace in a transversal $T$ the elements $1, 2, \cdots, \kappa$ into $\kappa$ ones then we have a permutation matrix of order $\kappa$.

\begin{defi}
A latin square $L$ of order $\kappa$ is called \textbf{resolvable} if it is the sum of $\kappa$ transversals.
\end{defi}

\begin{theo}
Let $L_{1}, L_{2}, \cdots, L_{\kappa - 1}$ be a complete set of $MPLS$ of order $\kappa \geq 3$ derived from a projective plane $\pi$ of order $\kappa$. Then $L_{1}$ is resolvable and is in fact $(\kappa - 2)$-fold resolvable.
\end{theo}

Proof: Let $L_{1} = [l_{ij}^{(1)}]$ and $L_{2} = [l_{ij}^{(2)}]$ be two mutually projective latin squares of order $\kappa$. The first row of $L_{2}$ has with each row of $L_{1}$ an element in common. Of these $k$ independent intersection points we can form a transversal $T_{1}$ of $L_{1}$ for these $\kappa$ intersection points occur in $\kappa$ different rows of $L_{1}$ and also in $\kappa$ columns of the first row of $L_{2}$. By repeating the same procedure with the second row of $L_{2}$ and its intersection points with the $\kappa$ rows of $L_{1}$ we find the second transversal $T_{2}$ of $L_{1}$, etc, so that we find that $L_{1}$ is resolvable. Because $L_{1}$ is also resolvable by the rows of the latin squares $L_{3}$,..., $L_{\kappa - 1}$ we see that $L_{1}$ is in fact $(\kappa - 2)$-fold resolvable. 

\section{More Food for the Fineproofer}

\newtheorem{theo3}[defi]{Theorem D of K\"{o}nig [1916]}
\begin{theo3}
If in a matrix $F = [f_{ij}]$ of order $n$, whose entries are nonnegative integers, every row and every column has the same positive sum, then at least one term of $det\ F$ is unequal to zero. 
\end{theo3}

\newtheorem{theo4}[defi]{Theorem E of K\"{o}nig [1916]}
\begin{theo4}
If in a matrix $F =[f_{ij}]$ of order $n$, every row and every column has the same positive number of $k$ nonzero entries, then $det\ F$ contains at least $k$ terms which are unequal to zero.
\newline (These $k$ terms can be chosen in such a way that every nonzero entry in the matrix $F$ appears in exactly one of the $k$ terms.) 
\end{theo4}

Frobenius[1917] was not happy with the graph-theoretic proof of Theorem D of K\"{o}nig and produced the following theorem: [but we know nowadays: \textit{In de wiskunde is alles geoorloofd als het maar klopt.}]

\newtheorem{theo5}[defi]{Theorem II of Frobenius [1917]}
\begin{theo5}
Let $F$ be a matrix of order $n$ and let $B = O$ be a submatrix of order $a \times b$ whose entries are zeroes and let $a + b = n + 1$. Then all terms of $det\ F$ are equal to zero.
\newline If all terms of $det\ F$ are equal to zero then there exists a submatrix $B = O$ of order $a \times b$ whose entries are all equal to zero and $a + b = n + 1$.
\end{theo5}

\newtheorem{defi3}[defi]{Definitions}
\begin{defi3}
A \textbf{binary} or \textbf{incidence matrix} $F$ of order $m \times n$ is an matrix, square or not, whose entries are zeroes or ones. A matrix with solely zeroes is denoted by $O$ and a matrix with solely ones is denoted by $J$. Some ones of an binary matrix $F$ are called \textbf{independent} if no two ones lie in the same row or in the same column. Let $v(F)$ denote the maximum number of independent ones. Further let a submatrix $B = O$ be of order $a \times b$ whose entries are solely zeroes with the sum $w(F) = w(B) = a + b$ being the greatest number for the matrix $F$ possible. If $F$ has no zeroes we put $w(F) = 0$. If we permute the rows and columns of matrix $F$ then the numbers $v(F)$ and $w(F)$ remain unchanged for the new matrix which we will call also $F$. We note that $v(F) \leq min (m, n)$.
\end{defi3}

\begin{theo}
Let $F = [f_{ij}]$ be an incidence matrix of order $m \times n$ whose entries are zeroes or ones. Suppose that after permutations of the rows and columns we can partition $F$ as follows: $F = \left[ \begin{array}{ll} A & B \end{array} \right]$ where the submatrix $A$ is of order $m \times (n - b)$ and the submatrix $B = O$ is of order $m \times b$ and the entries of $B$ are all zero. Then we have:
\newline $(a)$ If $w(F) = w(B) = m + b$ then $v(F) = v(A) = min(m, n - b)$.
\newline $(b)$ (Frobenius $[1917]$) If $w(A) > max (m, n - b)$ then $w(F) > w(B) = m + b$.
\newline $(c)$ If $w(F) = w(B) = m + b$ then $w(A) \leq  max (m, n - b)$.
\end{theo}

Proof: $(a)$ Let $v(A) < min (m, n - b)$ and $w(F) = w(B) = m + b$ then it follows from theorem 43 that $w(A) > max(m, n - b)$ but any zero matrix of $A$ can be extended by some rows of $b$ columns of zeroes of $B$ so $w(F) \geq w(A) + b > m + b$ or $w(A) + b > n - b + b \geq m + b$. In both cases we have a contradiction. This proves theorem $39(a)$.
$(b)$ If $w(A) > max (m, n - b)$ then $A$ contains after permutation a submatrix, let us say $E = O$, of order $e_{1} \times e_{2}$ whose entries are all zero and such that $w(E) = w(A) = e_{1} + e_{2}$. But the $e_{1}$ rows of $E = O$ meet $b$ columns of $B = O$. So $F$ contains a submatrix $E_{1} = O$ whose entries are all zero and is of order $e_{1} \times (e_{2} + b)$. Thus $w(F) \geq (e_{1} +(e_{2} + b)) = (e_{1} + e_{2}) + b = w(A) + b > max (m, n - b) + b$. Thus if $max (m, n - b) = m$ then $w(F) > m + b$ as stated or if $max (m, n - b) = n - b$ then $n - b \geq m$ and $w(F) > n \geq m + b$ as stated. $(c)$ is the contraposition of $(b)$.

\begin{theo}
If $P$ is a permutation matrix of order $n$ then $v(P) = n$ and $w(P) \leq n$.
\end{theo}

Proof: The $n$ ones are by definition independent, thus $v(P) = n$. We partition, after permutation, the matrix as $P = \left[ \begin{array}{ll} A & B \\ C & D \end{array} \right]$ where the submatrix $B$ is of order $a \times b$ and all elements of $B = O$ are zeroes and $w(P) = w(B)$ is a maximum number. Now in the rows of submatrix $A$ there are $a$ ones and in the columns of submatrix $D$ there are $b$ ones and $P$ contains $n$ ones. Thus $a + b \leq n$ and $w(P) = w(B) = a + b \leq n$.

\begin{theo}
Let $F$ be a binary matrix of order $n$. Then $v(F) = n \Leftrightarrow w(F) \leq n$.
\end{theo}

Proof: $(\Rightarrow)$ If $v(F) = n$ then $F$ can be written as the sum of a permutation matrix $P$ with $n$ independent ones and a binary matrix $M$: $F = P + M$. Let us say that $B = O$ is a submatrix of $F$ with solely zeroes and is of order $a \times b$ and such that $w(F) = w(B) = a + b$ is maximum. Now $P$ and $M$ contains a copy of $B = O$ but we know from theorem $40$ that $a + b \leq n$. Thus $w(F) \leq n$.
\newline \indent To prove $(\Leftarrow)$ we use induction on the number $n$. The theorem is true for $n = 1$, so we assume $n \geq 2$ and that matrix $F$ contains at least one one and one zero. We assume that $(i)$: $w(F) < n$ and $(ii)$: $w(F) = n$. In case of $(i)$ we choose an element one and its complement $G$ in matrix $F$. The submatrix $G$ is of order $n - 1$ and $w(G) \leq (n - 1)$ (otherwise $w(F) = n$). So by induction $v(G) = n - 1$ and $v(F) = 1 + (n - 1) = n$ and the theorem is proved. In case of $(ii)$ we can partition the matrix $F$ as follows: $F = \left[ \begin{array}{ll} A & B \\ C & D \end{array} \right]$ where the submatrix $B = O$ contains solely zeroes and is of order $x \times y$ and such that $w(F) = w(B) = x + y = n$. Then submatrix $A$ is of order $x$ and submatrix $D$ is of order $y$. Now $w(A) \leq x$ (otherwise $w(F) > n$ according to theorem $39 (b)$ which is a contradiction.) By induction we find $v(A) = x$. In the same manner we prove $v(D) = y$. Thus $v(F) = v(A) + v(D) = x + y = n$ and the theorem is proved. The theorem is equivalent with the contraposition of Theorem $II$ of Frobenius $[1917]$ $37$.

\begin{theo}
Let $F = [f_{ij}]$ be a binary matrix of order $m \times n$ whose entries are zeroes or ones. Then $v(F) = min (m, n) \Leftrightarrow w(F) \leq max (m, n)$.
\end{theo}

Proof: We assume $(i)$: $m < n$, $(ii)$: $m = n$ and $(iii)$: $m > n$. In case of $(i)$ we adjoin to matrix $F$ $(n - m)$ rows of solely ones to form a square matrix $H_{r}$ with the properties $v(H_{r}) = v(F) + (n - m) = m + (n - m) = n$ and $w(H_{r}) = w(F)$. Thus $v(F) = min (m, n) = m \Leftrightarrow v(H_{r}) = n \Leftrightarrow w(H_{r}) \leq n$ (according to theorem $41$) $\Leftrightarrow w(F) \leq n = max(n, m)$. This proves case $(i)$. The cases $(ii)$ and $(iii)$ are left to the reader.  

\begin{theo}
Let $F = [f_{ij}]$ be a binary matrix of order $m \times n$ whose entries are zeroes or
ones. Then $v(F) < min (m, n) \Leftrightarrow w(F) > max (m, n)$.
\end{theo}

Proof: This theorem is the contraposition of theorem 42.

\begin{theo}
Let $M$ be the incidence matrix of a finite projective plane $\pi$ of order $\kappa$. Then $M$ can be written as a sum of $\kappa + 1$ permutation matrices of order $\kappa^{2} + \kappa + 1$.
\end{theo}

Proof: Recall that in each row and each column of $M$ there are $\kappa + 1$ ones. We call the order of $M$ to be $n = \kappa^{2} + \kappa + 1$. After permutions of $M$ we can partitioned $M$ into the form $M = \left[ \begin{array}{ll} A & B \\ C & D \end{array} \right]$ where $B = O$ is a maximal submatrix whose entries are zeroes. The order of submatrix $A$ is $a \times (n - b)$, the order of submatrix $B$ is $a \times b$ and $w(M) = w(B) = a + b$, the order of submatrix $C$ is $(n - a) \times (n - b)$ and the order of submatrix $D$ is $(n - a) \times b$. We will prove that $w(M) = w(B) = a + b \leq n$. The number of ones in submatrix $A$ amount to $a (k + 1)$, the number of ones in submatrix $D$ amount to $b (k + 1)$ so that for submatrix $C$ there remains $n (k + 1) - (a + b) (k + 1)$ ones. But then $a + b \leq n$. Thus $w(M) = a + b \leq n$ and by theorem 41 we see that $v(M) = n$. The matrix $M$ contains $n$ independent ones or to put it in an different way $M$ can be written as $M = P_{1} + M_{1}$ where $P_{1}$ is an permutation matrix and $M_{1}$ is a binary matrix which have in each row and in each column $\kappa$ ones. In the same manner we can write $M = P_{1} + (P_{2} + M_{2})$ et cetera. The theorem can also be proved by Theorem $E$ of K\"{o}nig $[1916]$ $36$.

\begin{theo}
Let $L$ be a latin square of order $n$. Then any submatrix $A \subset L$ of order $a \times b$ with $a + b = n + m, 1 \leq m \leq n$, contains all elements $1, 2, \cdots, n$ at least $m$ times. 
\end{theo}

Proof: For $m = 1$ see L. Bleijenga[2002b]. The theorem is true for $a = n$ or $b = n$, so we assume $n > a, b > m \geq 1$. After permutation of $L$ we can assume that $L = \left[ \begin{array}{ll} A & B \\ C & D \end{array} \right]$. If a symbol $y$ occurs $m_{1}$ times in $A$, $0 \leq m_{1} < m$ then $y$ occurs in $a - m_{1}$ rows of $B$ and $b - m_{1}$ columns of $C$. Then $y$ occurs $n - a - b + m_{1} < 0$ times in $D$ which is impossible.

\begin{theo}
Let $G$ be a finite group of order $n$ and let $A$ and $B$ be subsets of $G$. If $|A| + |B| = n + 1$ then $AB = G$.
\end{theo} 

Proof: The multiplication table of the group $G$ is a latin square, let us say $L$, and $AB$ is a submatrix of $L$ whose order is $|A| \times |B|$ fulfills the condition $|A| + |B| = n + 1$ and according to theorem 45 for the case $m = 1$ the submatrix $AB$ contains all the elements of $G$. Thus $AB = G$.

\noindent \textbf{References}

\begin{description}
  \item[[ 1916]] D\'{e}nes K\"{o}nig, \"{U}ber Graphen und ihre Anwendung auf Determinantentheorie und Mengenlehre, Math. Ann. \textbf{77}, 453-465.
  \item[[ 1917]] Ferdinand Georg Frobenius, \"{U}ber zer\-leg\-bare Deter\-minan\-ten, Sitzungs\-berichte der K\"{o}nig\-lich Preu\ss ischen Aka\-demie der Wissen\-schaf\-ten zu Berlin 274-277 = FERDI\-NAND GEORG FRO\-BENIUS, GESAM\-MELTE AB\-HAN\-DE\-LUNGEN, III, Her\-aus\-gege\-ben von J-P. Serre, 701-704, 1968, \-Sprin\-ger-Verlag.
  \item[[ 1948]] N.G. de Bruijn and P. Erd\"{o}s, \textit{On a combinatorial problem}, Kon. Ne\-der\-l. Akad. Wetensch. Proc. A, \textbf{51}, (1948), 1277-1279 (=Indag. Math, \textbf{10}, (1948), 421-423).
  \item[[ 1953]] Lowell J. Paige and Charles Wexler, \textit{A Canonical Form for Incidence Matrices of Finite Projective Planes and their Associated Squares}, Portugaliae Mathe\-ma\-tica, Volume \textbf{12}, 1953, pp. 105-112.
  \item[[ 1993]] Problem Section, Problem 892, Proposed (and provided with solution) by L. Bleijenga, Nieuw Archief voor Wiskunde, Vierde serie, deel 11, \textbf{No. 3}, November 1993, p 284.
  \item[[ 1995]] Problem Section, Solutions, Problem 892, Solution by R.H. Jeurissen, Nieuw Archief voor Wiskunde, Vierde serie, deel 13, \textbf{No. 2}, July 1995, 244-245. 
  \item[[ 2002a]] Problem Section, Problem 27, Proposed by L. Bleijenga, Nieuw Archief voor Wiskunde, vijfde serie, deel 3, \textbf{nummer 1}, maart 2002, p. 96.
  \item[[ 2002b]] Problem Section, Solutions, Problem 27, Solution by L. Bleijenga, Nieuw Archief voor Wiskunde, vijfde serie, deel 3, \textbf{nummer 3}, september 2002, p. 280.
  \item[[ 2006]] Dhananjay P. Mehendale, Finite Projective Planes, arXiv:math.GM/\-0611492.
\end{description}
\end{document}